\documentclass[a4,12pt,multicol,makeidx]{article}
\def\origin{
  \clearpage
\vskip-\baselineskip\vskip-\topskip%
  \vbox to 0pt{\vskip-1in%
    \hbox to 0pt{\hskip-1in%
      \hbox to 0pt{\vrule width 1cm height .4pt depth 0mm\hss}%
      \vbox to 0pt{\hrule width .4pt height 0pt depth 1cm\vss}%
    \hss}%
  \vss}
  \vskip-\baselineskip
  \vbox to 0pt{\vskip-1in\vskip3cm%
    \hbox to 0pt{\hskip-1in\hskip3cm%
      \hbox to 0pt{\hss\vrule width 2cm height .4pt depth 0mm\hss}%
      \vbox to 0pt{\vss\hrule width .4pt height 1cm depth 1cm\vss}%
    \hss}%
  \vss}%
\vskip5mm\hskip10mm (3cm,3cm)
}%
   
\def\a{\alpha}   

\def\ou{\overline{u}} \def\uu{\underline{u}}    
\def\ov{\overline{v}}   
\def\hx{\hat{x}} \def\hy{\hat{y}} \def\hv{\hat{v}}

\def\tv{\tilde{v}} \def\ov{\overline{v}}
\def\ov{\overline{v}}   
 \def\l{\lambda}   \def\p{\partial}  \def\e{\varepsilon} 
\def\n{\nabla}    
\def\leq{\underline{<}} 
  \def\hx{\hat{x}} \def\hy{\hat{y}}
\def\n{{\nabla}}

\newenvironment{theorem}{%
\par \bigskip \it}{%
\bigskip \par}

\makeindex
\title{Strong maximum principle for radiative transfer type operators
}
\author{Mariko Arisawa\\ DAMTP, Centre for Mathematical Sciences
\\University of Cambridge\\
Wilberforce road, Cambridge\\
CB3 0WA England\\
E-mail: M.Arisawa@damtp.cam.ac.uk
}
\date{}
\begin{document}
\maketitle
\bigskip

{\bf Synopsis: } The strong maximum principle ((SMP) in short) for subsolutions of the radiative transfer type equations is shown in this paper. 
We treat a general class of integro-differential equations, defined in the product space of the space variable "$x$" and the velocity variable "$v$". 
 The equations consist of two terms : a nonlocal integral operator in $v$ variable, and a first-order partial-differential operator in $x$ variable. 
The nonlocal term represents the jump process in $v$ direction, and the term of the first-order partial differential operator describes the drift in $x$ 
direction. In particular, the drift in $x$ is generated by the velocity variable $v$. Based on the idea of the propagation of maxima, we give a general 
sufficient condition ((A) in the paper) so that the (SMP) holds for the present class of nonlocal equations. The framework of the viscosity solution 
is used to formulate the problem, and the related existence and uniqueness of solutions are also given. \\
\section{Introduction.} 
$\qquad$  The strong maximum principle concerning with the radiative transfer operator 
\begin{equation}\label{rt}
	-\left\langle v, \n_x u(x,v) \right\rangle - 
	\int_{V} [u(x,v') - u(x,v)] dq(v') \qquad (x,v)\in {\Omega} \times V
\end{equation}
is shown in this paper. Here, $\Omega$, $V$ are open domains in ${\bf R^N}$, $V$ is  precompact, and $dq(v')$ is a positive bounded Radon measure 
satisfying 
\begin{equation}\label{mass}
	\int_{V} 1 dq(v')=m>0 \quad (\hbox{constant}). 
\end{equation}
We generalize (\ref{rt}) to the following class of nonlocal integro-differential equations 
$$
	\sup_{\a \in \mathcal{A}} \{ -\left\langle b(x,v,\a), \n_x u(x,v) \right\rangle \} - \rho^{\gamma}
	\int_{V} [u(x,v') - u(x,v)] dq(v') =0 
	\qquad\qquad\qquad
$$
\begin{equation}\label{general}
	\qquad\qquad\qquad\qquad\qquad\qquad\qquad\qquad\qquad
	\qquad (x,v)\in {\Omega} \times V, 
\end{equation}
 where $\rho=\int_{V} |u(x,v')| dq(v')$, $\gamma\in {\bf R}$, $\mathcal{A}$ is a subset of a metric space, $b$ is defined in  ${\bf R^N} \times V\times {\mathcal A}$ with values in ${\bf R^N}$ satisfying 
\begin{equation}\label{b}
	|b(x,v,\a)-b(x',v',\a)|\leq C(|x-x'|+ |v-v'|)
\end{equation}
$$
	\qquad\qquad\qquad\qquad\qquad\qquad\qquad
	\quad \forall (x,v), (x',v')\in {\bf R^N} \times V, \forall \a \in {\mathcal A}. 
$$
Remark that (\ref{general}) contains the operator (\ref{rt}) as a special case, namely for $\mathcal{A}=\emptyset$, $b(x,v,\a)=v$ in ${\Omega} \times V$, and $\gamma=0$. We treat (\ref{general}) in the framework of viscosity solutions. We use the term the "strong maximum principle" ((SMP), in short) in the following sense. 
 Let $V_0=\hbox{supp}(dq(v'))$. \\

(SMP)  The equation (\ref{general}) is said to satisfy (SMP), if 
 for any subsolution $u$ of (\ref{general}) which  takes a maximum at $(x_0,v_0)$$\in \Omega\times V_0$,  $u(x,v)\equiv u(x_0,v_0)$ holds in $\Omega \times V_0$. \\

We establish (SMP) for (\ref{general}) in Theorem 2.1, under the condition (A) given in \S2 below. 
The strong maximum principle is known to hold, in the classical sense, for the second-order uniformly elliptic operator (Gilbarg and Trudinger \cite{gt}). And, it is also known to hold for the possibly degenerate elliptic operators in the framework of viscosity solutions, by Trudinger \cite{trudinger}, 
 Kawhol, and Kutev \cite{kk}, and Bardi and Da Lio \cite{bardi1}, \cite{bardi2}.  We also refer the readers to Bony \cite{bony} which made clear the relationship between (SMP), the propagation of maxima, and the hypoellipticity (in the sense of H{\"{o}rmander) of the second-order degenerate elliptic operator. Recently,  the author showed in \cite{loc}, \cite{qphomo}, that (SMP) holds for a class of integro-differential equations with L{\'e}vy operators  
\begin{equation}\label{levy}
	H(x,\n u,\n^2 u)- \int_{{\bf R^N}} [u(x+z)-u(x)-{\bf 1_{|z|\leq 1}}\left\langle \n u(x),z \right\rangle ] d\overline{q}(z) =0 \quad \hbox{in } \quad{\bf R^N},
\end{equation}
where $H$ is a fully nonlinear possibly degenerate second-order elliptic operator defined in $\Omega\times {\bf R^N}\times {\bf S^N}$ such that 
$$
	H(x,0,O)\geq 0 \qquad \forall x\in {\bf R^N}, 
$$
$d\overline{q}(z)$ is a positive Radon measure such that 
$$
	\int_{{\bf R^N}} \min\{ 1, |z|^2 \} d\overline{q}(z) <  \infty, 
$$
and that there exists a ball $B(0,r)\subset {\bf R^N}$, centered at the origin with radius $r>0$ such that 
\begin{equation}\label{trans}
	 B(0,r) \subset \hbox{supp}(d\overline{q}(z)). 
\end{equation}
Rmark that (\ref{general}) and (\ref{levy}) are completely different (compare also with another problem (\ref{rtlevy}) in below). 
The strong maximum principle for the radiative transfer equation can be applied to solve the ergodic problem (see Arisawa and Lions \cite{al}) and  the homogenization.  The homogenization of 
 the radiative transfer equation has been studied by Goudon and Poupaud \cite{poupaud},  Lions and Toscani \cite{lionstoscani}, and Bardos, Golse, and Perthame \cite{bgp},  to approximate the diffusion equation or the porous medium equation. 
 We  also refer to Evans \cite{ev1},  \cite{ev2} the treatment of the similar problems in the framework of viscosity solutions. Our motivation to study (SMP)  for (\ref{rt}) 
  is to treat the asymptotic problem in the presence of the controls, via the ergodicity. To be more precise (for the purpose of the application), the operators in \cite{bgp},  \cite{poupaud},  and \cite{lionstoscani}  are  in the following form. For $V_0=\hbox{supp} (dq(v'))$$\subset V$
\begin{equation}\label{rt0}
	-\left\langle v, \n_x u(x,v) \right\rangle - 
	\int_{V_0} [u(x,v') - u(x,v)] dq(v') \qquad (x,v)\in {\Omega} \times V_0. 
\end{equation}
The following are examples of $V_0$, the measure $dq(v)$, and the drift $b(x,v,\a)$ in the cited papers. \\

{\bf Example 1.1. (\cite{poupaud}, \cite{lionstoscani} )} $V_0=\{v^1,...,v^M\}$ ($M\in {\bf N}$ fixed, $v^i\in  {\bf R^N}$), $dq(v)=\sum_{i=1}^{M} w_i v^i$, where 
$\sum_{i=1}^{M} w_i =1$ ($w_i\geq 0$), and $b(x,v,\a)=v$. In this case, (\ref{general}) becomes the system of $M$ equations for $u_i(x)=u(x,v^i)$:
$$
	-\left\langle v^i, \n_x u_i\right\rangle - (\sum_{j=1}^M w_ju_j -u_i)=0 \qquad \hbox{for}\quad 1\leq i\leq M.
$$

{\bf Example 1.2. (\cite{bgp}, \cite{poupaud})} $V_0={\bf S^{N-1}}$, $dq(v)=dv$ (the normalized Lebesgue measure on the sphere), and $b(x,v,\a)=v$.\\

Remark that the extention of the domain $V_0$ to $V$, the change of the operator from (\ref{rt0}) to (\ref{rt}), and the reformulation of the problem from $\Omega \times V_0$ to $\Omega \times V$ do not influence the radiative transfer model, for the subset $V_0$$=\hbox{supp} (dq(v'))$$\subset V$ is invariant under the jump 
process generated by $ -\int_{V} [u(x,v') - u(x,v)] dq(v')$.  
To contrast with the operator $\int_{V}[u(x,v')-u(x,v)]dq(v')$,  for a special case of $b(x,v,\a)=v$,  we  also treat 
$$
	-\left\langle v, \n_x u(x,v) \right\rangle  - \int_{{\bf R^N}} [u(x,v+w) - u(x,v)] dq(w) =0 
	\qquad\qquad\qquad
$$
\begin{equation}\label{rtlevy}
	\qquad\qquad\qquad\qquad\qquad\qquad\qquad\qquad\qquad
	\qquad (x,v)\in {\Omega} \times {\bf R^N}, 
\end{equation}
where the non-local term is the L{\'e}vy type operator in $v$, and $dq(\cdot)$ is a positive bounded Radon measure which satisfies (\ref{mass}). Differently from 
 (\ref{rt}), for the operator $\int_{{\bf R^N}} [u(x,v+w) - u(x,v)] dq(w)$, $W_0=\hbox{supp} (dq(w))$  is the set of vectors of jumps, and is 
not the invariant set of the corresponding  jump process in the phase space. We show (SMP) for (\ref{rtlevy}) in Theorem 2.5 below. \\ 

We denote $USC(\Omega\times V)$ (resp. $LSC(\Omega\times V)$) the set of upper (resp. lower) semicontinuous functions in $\Omega\times V$.
The set of all second-order subdifferentials (resp. superdifferentials) of $u\in USC(\Omega\times V)$ (resp. $w\in LSC(\Omega\times V)$) at $(x,v)\in \Omega\times V$ is denoted by  $J^{2,+}_{\Omega\times V} u(x,v)$ (resp. $J^{2,-}_{\Omega\times V} w(x,v)$). We refer the readers to Crandall, Ishii, and Lions \cite{users} for the elementary definitions and notations in the theory of the viscosity solution. The following is the definition of the viscosity solution for (\ref{general}). \\

{\bf Definition 1.1.} A function $u\in USC(\Omega\times V)$ (resp. $w\in LSC(\Omega\times V)$) is said to be a viscosity subsolution (resp. supersolution) of (\ref{general}) 
  if  for any $(p,X)\in J^{2,+}_{\Omega\times V}u(x,v)$ (resp. $(p,X)\in J^{2,-}_{\Omega\times V}w(x,v)$), the following holds. 
$$
	\sup_{\a \in \mathcal{A}} \{ -\left\langle b(x,v,\a), p \right\rangle \} - 
	\rho^{\gamma}\int_{V} [u(x,v') - u(x,v)] dq(v') \leq 0. 
$$
(resp. 
$$
	\sup_{\a \in \mathcal{A}} \{ -\left\langle b(x,v,\a), p \right\rangle \} - 
	\rho^{\gamma}\int_{V} [u(x,v') - u(x,v)] dq(v') \geq 0.)
$$
If $u$ is a viscosity subsolution and a viscosity supersolution at the same time, $u$ is said to be a viscosity solution.\\

It is not so difficult to verify the comparison and the existence of the viscosity solutions for (\ref{general}) with $\gamma=0$  with some  boundary conditions (Dirichlet, periodic, etc). We state some results with the Dirichlet boundary condition in \S3  Appendix,  to justify our (SMP) results in the framework of viscosity solutions. For (\ref{rtlevy}), we refer the readers to  Alvarez and Tourin \cite{at}, Arisawa \cite{new}, \cite{cor}, \cite{def}, Barles, Buckdahn, and Pardoux \cite{bbp}, Barles and Imbert \cite{bi}. For  (\ref{general}) with $\gamma\neq 0$, the existence and the uniqueness of solutions are known in \cite{lionstoscani} and  \cite{bgp}, by different approachs.  Here, we do not enter in details for the comparison and the existence of the solutions for such a case, and we concentrate on the study of (SMP).\\

\section{Strong maximum principle}

	$\qquad$ Consider the following deterministic system in ${\bf R^N}$. For fixed measurable vector valued functions $\a(\cdot)$ : 
$[0,\infty)$ $\to {\mathcal A}$ and $\tilde{v}(\cdot)$$:[0,\infty)\to V_0$, 
\begin{equation}\label{system}
	dX_{\a,\tv}(t)=b(X_{\a,\tv}(t),\tv(t),\a(t)) dt \quad t>0; \qquad X_{\a,\tv}(0)=x. 
\end{equation}
Remark that for any fixed $\a(\cdot)$$:[0,\infty)\to \mathcal{A}$ and $\tilde{v}(\cdot)$$:[0,\infty)\to V_0$ measurable, $\{X_{\a,\tv}(t)\}_{t\geq 0}$ is a trajectory in ${\bf R^N}$ and not in ${\bf R^N}\times V_0$ (or in ${\bf R^N}\times V$). \\

(A) The deterministic controlled system (\ref{system}) is said to be controllable in $\Omega$ if for any $x$, $y\in \Omega$ there exist $\tv(\cdot)$$:[0,\infty)\to V_0$, and $\a(\cdot)$$:(0,\infty)$$\to {\mathcal A}$,  and $T_{x,y}\geq 0$ such that the solution $X_{\a,\tv}(t)$ of (\ref{system}) satisfies $X_{\a,\tv}(0)=x$ and  
$X_{\a,\tv}(T_{x,y})=y$. \\

If $b(x,v,\a)$ does not depend on $v$ or $\a$, we omit the index and write $X_{\a}(t)$ or $X_{v}(t)$.\\

{\bf Example 2.1.}  Let $N=1$, $\Omega=(-1,1)$, $V=(-2,2)$, $V_0$$=\{-1,1\}$, ${\mathcal A}=\emptyset$, $b(x,v,\a)$$=v$ for any $x\in \Omega$, 
$dq(v')=\frac{1}{2}(\delta_{-1}+\delta_{1})$. Then,  for any 
$x$, $y$$\in (-1,1)$ we can take either $\tv\equiv -1$ or $\tv\equiv 1$ ($\forall t\geq 0$) so that the solution $X_{\tv}(t)=x\pm t$ ($t\geq 0$) of (\ref{system}) 
 satisfies (A), i.e. $X_{\tv} (0)=x$, $X_{\tv} (T_{xy})=y$ for $T_{xy}=|y-x|$. \\

{\bf Example 2.2.} (Communicated by P.-L. Lions.) Let $N=2$, $\Omega=(0,1)^2$, $V=B(0,2)$, $V_0$$=\{\pm{\bf e_1}, \pm{\bf e_2}\}$,  
where ${\bf e_1}=(1,0)$, ${\bf e_2}=(0,1)$, 
 ${\mathcal A}$$=\emptyset$,  $b(x,v,\a)$$=v$ for any $(x,v)\in \Omega\times V$, 
 and 
 $dq(v')$$=\frac{1}{4}(\delta_{{\bf e_1}}$$+\delta_{-{\bf e_1}}$$+ \delta_{{\bf e_2}}+$$\delta_{-{\bf e_2}})$. Then, for any $x=(x_1,x_2)$, $y=(y_1,y_2)$$\in \Omega$, we can take 
\begin{eqnarray}
	\tv(t) &=& \frac{y_1-x_1}{|y_1-x_1|} {\bf e_1} \quad 0\leq t\leq |y_1-x_1| \nonumber \\
	   \quad    &=& \frac{y_2-x_2}{|y_2-x_2|}  {\bf e_2} \quad   |y_1-x_1|\leq t \leq \sum_{i=1}^2|y_i-x_i|. \nonumber 
\end{eqnarray}
so that  the solution $X_{\tv}(t)$  of (\ref{system}) satisfies (A), i.e. $X_{\tv} (0)=x$, $X_{\tv} (T_{xy})=y$ for $T_{xy}=\sum_{i=1}^2|y_i-x_i|$. \\

{\bf Example 2.3.} Let $N\geq 1$, $\Omega=B(0,1)$, $V=B(0,2)$, $V_0$$=B(0,1)$,
 ${\mathcal A}$$=\{\a\in {\bf R^N}|\quad  |\a|=1,\quad \hbox{or}\quad \a=0\}$,  $b(x,v,\a)$$=\a$ for any $(x,v)\in \Omega\times V$, 
 and 
 $dq(v')=dv'$ (the normalized Lebesgue measure). Then, for any $x$, $y$$\in \Omega$ ($x\neq y$), by taking 
$\a=\frac{y-x}{|y-x|}$, $T_{xy}=|y-x|$, $X_{\a}(t)=x+\a t$  ($t\in [0,T]$) satisfies (\ref{system}) with $X_{\a}(0)=x$, $X_{\a}(T_{xy})=y$. \\

	The following is our main result of (SMP) for (\ref{general}). \\

{\bf Theorem 2.1.$\quad$}
\begin{theorem}  Let (\ref{mass}) and (\ref{b}) hold in (\ref{general}). Assume that the controlled deterministic system (\ref{system}) satisfies the condition (A), i.e. it is controllable in $\Omega$. 
Let $u\in USC(\Omega\times V)$ be a subsolution of   (\ref{general}). Assume that $u$ takes a maximum at 
$(x_0,v_0)$$\in \Omega\times V_0$, i.e. $u(x_0,v_0)=$$\max_{\Omega\times V_0}u(x,v)$. Then, $u(x,v)=u(x_0,v_0)$ for any $(x,v)$$\in \Omega\times V_0$. \\
\end{theorem}

{\bf Proof.} We may assume that $m=1$ in (\ref{mass}), for otherwise we devide the both hand sides by 
 $m>0$. We prove the claim  in the following three steps.\\
 (Step 1.)  Denote $Y$ the set of points $y\in \Omega$ which is controllable from $x_0$, i.e. 
$$
	Y=\{y\in \Omega|\quad \exists \tv(\cdot): [0,\infty)\to V_0,\quad \exists \a(\cdot): [0,\infty)\to {\mathcal A},\quad \exists T_{xy}\geq0
$$
$$
	\hbox{s.t.}\quad 
    X_{\a,\tv}(t):\quad dX_{\a,\tv}(t)=b(X_{\a,\tv}(t),\tv(t),\a(t)) dt \quad t>0; \qquad X_{\a,\tv}(0)=x
$$
$$
	\hbox{satisfies}\quad X_{\a,\tv}(0)=x_0,\quad X_{\a,\tv}(T_{xy})=y 
	\}. \qquad\qquad\qquad
$$
Let $Y_0=\{x_0\}$. 
Let $Y_1$ be the set of points $y\in \Omega$ which is controllable from $x_0$ by a control $\a(\cdot)$ and a constant control $\ov(t)\equiv v$ ($v\in V_0$, $t\geq 0$), i.e. 
$$
	Y_1=\{y\in \Omega|\quad \exists \ov(t)\equiv v\in V_0 \quad (\forall t\geq 0),\quad \exists \a(\cdot): [0,\infty)\to {\mathcal A},\quad \exists T_{xy}\geq0
$$
$$
	\hbox{s.t.}\quad 
    X_{\a,\ov}(t):\quad dX_{\a,\ov}(t)=b(X_{\a,\ov}(t),\ov(t),\a(t)) dt \quad t>0; \qquad X_{\a,\ov}(0)=x
$$
$$
	\hbox{satisfies}\quad X_{\a,\ov}(0)=x_0,\quad X_{\a,\ov}(T_{xy})=y 
	\}. \qquad\qquad\qquad
$$
Define $Y_k$ ($k\geq 2$) inductively the set of points $y\in \Omega$ which is controllable from points $x\in Y_{k-1}$ by a control $\a(\cdot)$ and a 
constant control $\ov(t)\equiv v$ ($v\in V_0$, $t\geq 0$), i.e. 
$$
	Y_k=\{y\in \Omega|\quad \exists \ov(t)\equiv v\in V_0 \quad (\forall t\geq 0),\quad \exists x\in Y_{k-1} \exists \a(\cdot): [0,\infty)\to {\mathcal A},\quad \exists T_{xy}\geq0
$$
$$
	\hbox{s.t.}\quad 
    X_{\a,\ov}(t):\quad dX_{\a,\ov}(t)=b(X_{\a,\ov}(t),\ov(t),\a(t)) dt \quad t>0; \qquad X_{\a,\ov}(0)=x
$$
$$
	\hbox{satisfies}\quad X_{\a,\ov}(0)=x,\quad X_{\a,\ov}(T_{xy})=y 
	\}. \qquad\qquad\qquad
$$
Clearly $Y_0\subset Y_1\subset Y_2\subset ... \subset Y_k \subset ...$. 
Since the measurable function $\tv(t)$: $[0,T]\to {V_0}$ ($T>0$ ) is the uniform limit of a sequence of the piecewise constant functions, say $\{\ov_n(t)\}_{n\in {\bf N}}$  ($\ov_n(t)\in V_0$, $\forall t\in [0,T]$), we know that 
\begin{equation}\label{Y}
	Y=\cup_{k\geq 0} Y_k. 
\end{equation}
(Step 2.) Let $M=u(x_0,v_0)=\max_{\Omega\times V} u(x,v)$. We shall see 
\begin{equation}\label{M}
	u(x_0,v)=M \quad \forall v\in {V_0}=supp(dq(v')).
\end{equation}
Since $(0,O)\in J^{2,+}_{\Omega}u(x_0,v_0)$, from Definition 1.1, 
$$
	\sup_{\a \in \mathcal{A}} \{ -\left\langle b(x_0,v_0,\a), 0 \right\rangle \} - \rho_0^{\a}
	\int_{V} [u(x_0,v') - u(x_0,v_0)] dq(v') \leq 0,
$$
where $\rho_0=\int_{V} |u(x_0,v')| dq(v')$. 
Assume that the set $V_1=\{v\in V_0| u(x_0,v)<M\}$ has a positive measure. Since $u(x_0,v)-u(x_0,v_0)\leq 0$ ($\forall v\in V$), 
 and since $\rho_0>0$, the above leads 
$$
	0< -\int_{V_1} u(x_0,v')-u(x_0,v_0) dq(v') \leq 
   -  \int_{V} u(x_0,v')-u(x_0,v_0) dq(v') \leq 0,
$$
which is a contradiction. Thus, (\ref{M}) was shown. \\

(Step 3.) We show the following by induction with respect to $k$. 
\begin{equation}\label{yk}
	u(x,v)=M \quad \forall (x,v)\in Y_k \times V_0, \quad \forall k=0,1,2,...
\end{equation}
In fact, from Step 2 (\ref{yk}) is true for $k=0$. By assuming that (\ref{yk}) holds for $k\in {\bf N}$, we see that it also holds for $k+1$. 
Put $g(x)=\int_{V} u(x,v') dq(v')$. We rewrite (\ref{general})  to
\begin{equation}\label{hj}
	u(x,v)+ \sup_{\a \in \mathcal{A}} \{ -\left\langle b(x,v,\a), \n_x u(x,v) \right\rangle \} - g(x)
	=0 \qquad (x,v)\in {\bf R^N} \times V. 
\end{equation}
Let $x_k\in Y_k$. From the assumption $u(x_k,v)=M$ for any $v\in V_0$. 
Let  $\ov(t)\equiv v\in V_0$, and let $\a(\cdot)$$:[0,\infty)\to {\mathcal A}$ be fixed temporarily, which we choose later.  Consider : 
$$
	dX_{\a,\ov}(t)=b(X_{\a,\ov}(t),v,\a(t)) dt \quad t>0; \qquad X_{\a,\ov}(0)=x_k. 
$$
Denote the exit time of $X_{\a,\ov}(t)$ from $\Omega$ as: 
$$
	\tau^{\a}_{x_k}=\min\{t\geq 0|\quad X_{\a,\ov}(t)\in \Omega^{c}\}. 
$$
Since (\ref{hj}) is 
the infinite horizon  Hamilton-Jacobi equation (see \cite{users}), we have 
$$
	u(x_k,v)=\inf_{\a(\cdot)} \{ \int_0^{\tau^{\a}_{x_k}} e^{-t} g(X_{\alpha,\ov}(t))dt + e^{-\tau^{\a}_{x_k}}u(X_{\alpha,\ov}(\tau^{\a}_{x_k}),v)\}
	\qquad\qquad\qquad\qquad\qquad
$$
$$\quad
	\leq  \int_0^{\tau^{\a}_{x_k}} \int_V e^{-t} u(X_{\alpha,\ov}(t),v') dtdq(v') + e^{-\tau^{\a}_{x_k}} M \quad \hbox{for}\quad \forall \a(\cdot):[0,\infty)\to {\mathcal A}
$$
$$
	\leq \int_0^{\tau^{\a}_{x_k}} e^{-t}M dt + e^{-\tau^{\a}_{x_k}} M =M. \qquad\qquad\qquad\qquad\qquad\qquad\qquad\qquad
$$
Since $u(x_k,v)=M$,  we get 
\begin{equation}\label{mM}
	0\leq \int_0^{\tau^{\a}_{x_k}} \int_V e^{-t} [u(X_{\alpha,\ov}(t),v')-M] dtdq(v') \leq 0  \quad \hbox{for}\quad \forall \a(\cdot):[0,\infty)\to {\mathcal A}. 
\end{equation}
Since $u(X_{\alpha,\ov}(t),v')-M\leq 0$ for any $v'\in V_0$, the above leads 
\begin{equation}\label{aconst}
	u(X_{\a,\ov}(t),v')=M \quad \hbox{a.e.}\quad t\in[0,\tau^{\a}_{x_k}],\quad v'\in V_0,\quad \forall \a(\cdot): [0,\infty)\to {\mathcal A}. 
\end{equation}
From the definition, for any $x_{k+1}\in$$Y_{k+1}$ there exists $x_k\in Y_k$ such that for some $\a(\cdot)$, $\ov$, and $T>0$ the solution $X_{\a,\ov}(t)$ 
 of (\ref{system}) satisfies $X_{\a,\ov}(0)=x_k$ and $X_{\a,\ov}(T)=x_{k+1}$.  
Therefore, (\ref{aconst}) leads 
$$
	u(x_{k+1},v)=M \quad \forall v\in V_0. 
$$
By induction, we thus proved (\ref{yk}). \\

The opposite sign of the above result does not hold in general : even  if $u\in LSC(\Omega\times V)$ is a supersolution of (\ref{general}) and $u$ takes a minimum at $(x_0,v_0)$$\in \Omega\times V_0$, $u(x,v)=u(x_0,v_0)$ ($\forall (x,v)\in \Omega\times V_0$) does not hold in general. Instead of (\ref{general}), if we consider 

$$
	\inf_{\a \in \mathcal{A}} \{ -\left\langle b(x,v,\a), \n_x u(x,v) \right\rangle \} - \rho^{\gamma}
	\int_{V} [u(x,v') - u(x,v)] dq(v') =0 
	\qquad\qquad\qquad
$$
\begin{equation}\label{general2}
	\qquad\qquad\qquad\qquad\qquad\qquad\qquad\qquad\qquad
	\qquad (x,v)\in {\Omega} \times V, 
\end{equation}

the following holds. \\

{\bf Proposition 2.2.$\quad$}
\begin{theorem}  Let (\ref{mass}) and (\ref{b}) hold in (\ref{general2}). Assume that the controlled deterministic system (\ref{system}) satisfies the condition (A), i.e. it is controllable in $\Omega$. 
Let $u\in LSC(\Omega\times V)$ be a supersolution of   (\ref{general2}). Assume that $u$ takes a minimum at 
$(x_0,v_0)$$\in \Omega\times V_0$,  i.e. $u(x_0,v_0)=$$\min_{\Omega\times V_0}u(x,v)$. Then, $u(x,v)=u(x_0,v_0)$ for any $(x,v)$$\in \Omega\times V_0$. \\
\end{theorem}

The proof of Proposition 2.2 is quite similar to that of Theorem 2.1, and we do not write it here. 
From Theorem 2.1, we have the following. To simplify the situation, we assume that $\Omega$ and $V$ are the tori. \\

{\bf Theorem 2.3.$\quad$}
\begin{theorem} Let $\Omega\times V={\bf T^{N}}\times {\bf T^{N}}$. 
Let (\ref{mass}) and (\ref{b}) hold in (\ref{general}). Assume that the controlled deterministic system (\ref{system}) satisfies the condition (A), i.e. it is controllable in $\Omega$. 
Let $u$ be a solution of  (\ref{general}). Assume that $u$ takes a maximum at 
$(x_0,v_0)$$\in \Omega\times V_0$,  i.e. $u(x_0,v_0)=$$\max_{\Omega\times V_0}u(x,v)$. Then, $u(x,v)=u(x_0,v_0)$ for any $(x,v)$$\in \Omega\times V$. \\
\end{theorem}

{\bf Proof.} As before we may assume that $m=1$. Let $u(x_0,v_0)=M$. 
From Theorem 2.1, since $u(x,v)=M$ ($\forall (x,v)\in \Omega\times V_0$), (\ref{general}) becomes 
$$
	M^{\a}u(x,v)+ \sup_{\a \in \mathcal{A}} \{ -\left\langle b(x,v,\a), \n_x u(x,v) \right\rangle \} - M^{\a+1} =0. 
$$
Since $\Omega\times V$ is the torous, from the formula of the value function of the Hamilton-Jacobi equation, we get  
\begin{equation}\label{expression}
	u(x,v)=\int_{0}^{\infty} e^{-M^{\a} t} M^{\a+1} dt =M \qquad \forall (x,v)\in \Omega\times V. 
\end{equation}
The claim is thus provesd.\\

{\bf Remark.}  The proof of Theorem 2.3 relies on the explicite formula (\ref{expression}). From this point of view, we can replace the periodic BC to the 
 Neumann type BC or the state constraint BC on $\p(\Omega\times V)$ to have the same claim. \\

{\bf Example 2.4.} Let $N>1$, $\Omega={\bf R^N}$, $V=B(0,R)$ ($R>1$), $V_0$$={\bf S^{N-1}}$, $dq(v')=dv'$ (the normalized Lebesgue measure on ${\bf S^{N-1}}$), 
 ${\mathcal A}$$=\emptyset$, and $b(x,v,\a)$$=v$ ($\forall (x,v)\in \Omega\times V$).  Let $u\in USC(\Omega\times V)$ be a subsolution of 
$$
	-\left\langle v, \n_x u(x,v) \right\rangle - 
	\int_{V} [u(x,v') - u(x,v)] dq(v') = 0 \qquad (x,v)\in {\bf R^N} \times V. 
$$
Then,  if $u$ takes a maximum at $(x_0,v_0)$$\in \Omega\times V_0$, $u$ is constant in $\Omega\times V_0$. It is easy to check that  the system (\ref{system}) 
satisfies the controllability (A), and the claim follows from Theorem 2.1.\\

The following is a counter example. \\

{\bf Example 2.5.}  Let $N=1$, $\Omega={\bf R}$, $V=V_0=(0,R)$ $(R>0)$, $b(x,v,\a)=v$ ($\forall (x,v)\in \Omega\times V$). Consider 
\begin{equation}\label{counter}
	-\left\langle v, \n_x u(x,v) \right\rangle - 
	\int_{V} [u(x,v') - u(x,v)] dq(v') =0 \qquad (x,v)\in {\bf R} \times (0,R). 
\end{equation}
Set $u(x,v)=1$ ($x\geq 0$); $u(x,v)=1+x$ ($x\leq 0$). Then, $u$ is a continuous subsolution of (\ref{counter}), which  takes a maximum at 
$(x,v)$ ($\forall x\geq 0$, $\forall v\geq 0$). However, the system (\ref{system}) does not satisfy the condition (A), and $u$ is not constant.\\

We end the analysis of (SMP) for (\ref{general}) with the following proposition.  \\

{\bf Proposition 2.4.$\quad$}
\begin{theorem} 
Let (\ref{mass}) and (\ref{b}) hold in (\ref{general}). 
Let $u\in USC(\Omega\times V)$ be a subsolution of   (\ref{general}). Let $Z_0=$$\{x \in \Omega|\quad \exists v\in V_0 \quad \hbox{s.t.}\quad u(x,v)=\max_{\Omega\times V_0}u(x,v)=M\}$.  Let 
$$
	Z=\{y\in \Omega|\quad \exists x\in Z_0, \exists \a(\cdot): (0,\infty)\to \mathcal{A},\quad \exists \tv(\cdot): (0,\infty)\to {V_0}, \exists T_{x,y}\geq 0, 
$$
$$
	\quad \hbox{s.t.}\quad X_{\a,\tv}(t): 
	\frac{dX_{\a,\tv}(t)}{dt}=b(X_{\a,\tv}(t),\tv(t),\a(t)) \quad t>0; \quad
	X_{\a,\tv}(0)=x \qquad
$$
$$
	\hbox{satisfies}\quad X_{\a,\tv}(T_{x,y})=y\}. \qquad\qquad\qquad\qquad
$$
Then, $u(x,v)=M$ for any $(x,v)\in \overline{Z}\times V_0$.
\end{theorem}
{\bf Proof.} As in Theorem 2.1, we define the increasing sequence of sebsets of $\Omega$: 
let $Z_1$ be the set of points $y\in \Omega$ which is controllable from points in  $Z_0$ by a control $\a(\cdot)$ and a constant control $\ov(t)\equiv v$ ($v\in V_0$, $t\geq 0$), i.e. 
$$
	Z_1=\{y\in \Omega|\exists x\in Z_0, \exists \ov(t)\equiv v\in V_0  (\forall t\geq 0),\exists \a(\cdot): [0,\infty)\to {\mathcal A},\exists T_{xy}\geq 0
$$
$$
	\hbox{s.t.}\quad 
    X_{\a,\ov}(t):\quad dX_{\a,\ov}(t)=b(X_{\a,\ov}(t),\ov(t),\a(t)) dt \quad t>0; \qquad X_{\a,\ov}(0)=x
$$
$$
	\hbox{satisfies}\quad X_{\a,\ov}(T_{xy})=y 
	\}. \qquad\qquad\qquad
$$
Define $Z_k$ ($k\geq 2$) inductively the set of points $y\in \Omega$ which is controllable from points $y\in Z_{k-1}$ by a control $\a(\cdot)$ and a constant control  $\ov(t)\equiv v$ ($v\in V_0$, $t\geq 0$), i.e. 
$$
	Z_k=\{y\in \Omega|\exists x\in Z_{k-1}, \exists \ov(t)\equiv v\in V_0 (\forall t\geq 0), \exists \a(\cdot): [0,\infty)\to {\mathcal A},\exists T_{xy}\geq0
$$
$$
	\hbox{s.t.}\quad 
    X_{\a,\ov}(t):\quad dX_{\a,\ov}(t)=b(X_{\a,\ov}(t),\ov(t),\a(t)) dt \quad t>0; \qquad X_{\a,\ov}(0)=x
$$
$$
	\hbox{satisfies}\quad  X_{\a,\ov}(T_{xy})=y 
	\}. \qquad\qquad\qquad
$$
As before $Z_0\subset Z_1\subset Z_2\subset ... \subset Z_k \subset ...$, and 
\begin{equation}\label{Z}
	Z=\cup_{k\geq 0} Z_k. 
\end{equation}
 By using the argument in Step 2 of Theorem 2.1, we can show that 
$$
	u(x,v)=M\quad  \forall (x,v) \in Z_0 \times V_0.
$$ 
Then, by the similar induction in  Step 3 of Theorem 2.1, 
$$
	u(x,v)=M \quad \forall (x,v)\in Z_k \times V_0, \quad k=0,1,2,.... 
$$
Therefore, from (\ref{Z}) and from the upper semicontinuity of $u$, we proved the claim. \\

{\bf Example 2.6.}  Let $N=2$,  $\Omega={\bf T^2}$,  $V=(-R,R)^2$ $(R>0)$, $\mathcal{A}=\emptyset$, $b(x,v,\a)=(1,\gamma)$ ($\gamma$ is a fixed irrational number) $\forall (x,v)\in \Omega\times V$,  
 $dq(v')=dv'$ (the normalized Lebesgue measure on $(-R,R)^2$).  
 Let $u\in USC(\Omega\times V)$ be a subsolution of 
$$
	-\left\langle b, \n_x u(x,v) \right\rangle - 
	\int_{V} [u(x,v') - u(x,v)] dq(v') = 0 \qquad (x,v)\in {\bf T^2} \times V. 
$$
Assume that $u$ takes a maximum at $(x_0,v_0)$$\in \Omega\times V$. Then, $u$ is constant in $\Omega\times V$. In fact,   the dynamical system $\frac{dX(t)}{dt}=(1,\gamma)$ ($t>0$) is ergodic in ${\bf T^2}$, and thus $\overline{Z}=\Omega$. Therefore, from Proposition 2.4 $u$ is constant in $\Omega\times V$. \\

If the nonlocal term is the L{\'e}vy type operator, the following (SMP) holds. \\

{\bf Theorem 2.5.$\quad$}
\begin{theorem}
 Let (\ref{mass}), (\ref{b}), and (\ref{trans}) hold in (\ref{rtlevy}). Assume that the controlled deterministic system (\ref{system}) satisfies the condition (A), i.e. it is controllable in $\Omega$. 
Let $u\in USC(\Omega\times {\bf R^N})$ be a subsolution of   (\ref{rtlevy}). Assume that $u$ takes a maximum at 
$(x_0,v_0)$$\in \Omega\times {\bf R^N}$,  i.e. $u(x_0,v_0)=$$\max_{\Omega\times V}u(x,v)$. Then, $u(x,v)=u(x_0,v_0)$ for any $(x,v)$$\in \Omega\times {\bf R^N}$. 
\end{theorem}

{\bf Proof.} (Step 1.)  Let $M=u(x_0,v_0)=\max_{\Omega \times  {\bf R^N}} u(x,v)$. Set $Z_0=\{x \in \Omega|\quad \exists v\in {\bf R^N} \quad \hbox{s.t.}\quad u(x,v)=M\}$. We show
\begin{equation}\label{last}
	u(x,v)=M \quad \forall (x,v)\in Z_0\times {\bf R^N}. 
\end{equation}
 We assume that the claim is not true, and shall get a contradiction.  So,   assume that  for some $x_1\in Z_0$,
$V_1=\{v\in {\bf R^N}|\quad u(x_1,v)=M\}\neq \emptyset$, and  $V_2=V_1^c \neq \emptyset$, too. For $r>0$ in (\ref{trans}), take $v_2\in V_2$ (i.e. $u(x_1,v_2)< M$),  such that $\hbox{dist}(v_2,V_1)< r$. 
Take $v^{\ast}\in V_1$ (i.e. $u(x_1,v^{\ast})=M$) such that $|v^{\ast}-v_2|<r$. Then, since $u(x_1,v^{\ast})=$$\max_{\Omega\times {\bf R^N}}u(x,v)$, $(0,O)\in J^{2,+}_{\Omega\times {\bf R^N}}$$u(x_1,v^{\ast})$,  and the definition of the viscosity subsolution leads 
$$
	0\leq 0-\int_{B(0,r)} u(x_1,v^{\ast}+w)-u(x_1,v^{\ast}) dq(w)
	\leq \int_{{\bf R^N}} u(x_1,v^{\ast}) - u(x_1,v^{\ast}+w)dq(w)
	\leq 0. 
$$
Thus, $u(x_1,v^{\ast}+w)=M$ almost evelywhere in $w\in B(0,r)$. However this contradicts to $u(x_1,v_2)<M$, for $|v^{\ast}-v_2|<r$ and for $u$ is upper semicontinuous. 
Therefore, (\ref{last}) was proved.\\
(Step 2). By using (\ref{last}), we repeat the argument in Theorem 2.1 to get $u(x,v)=M$ for any $(x,v)\in \Omega\times {\bf R^N}$.  \\

\section{Appendix : comparison principle and existence of viscosity solutions.} 

$\qquad$  In this section, we briefly show the comparison and the existence of the viscosity solutions for (\ref{general}) when $\gamma= 0$. 
We treat the stationary problem. The evolutionary problem can be treated similarly.
Let $\lambda>0$, and consider 
$$
	\l u + \sup_{\a \in \mathcal{A}} \{ -\left\langle b(x,v,\a), \n_x u(x,v) \right\rangle\} - 
	\int_{V} [u(x,v') - u(x,v)] dq(v') 
$$
\begin{equation}\label{stationary}
\qquad\qquad\qquad\qquad\qquad\qquad\qquad\qquad
	- g(x,v)=0 \qquad \hbox{in} \quad \Omega\times V,
\end{equation}
\begin{equation}\label{bc}
	u=\psi(x,v) \qquad \hbox{on} \quad \p(\Omega\times V),
\end{equation}
where $\int 1 dq(v')=1$,  and for some $\theta\in (0,1]$, $h=g$ or $\psi$ satisfies 
\begin{equation}\label{g}
	|h(x,v)-h(x',v')|\leq C(|x-x'|^{\theta}+|v-v'|^{\theta})\quad (x,v), (x',v')\in \Omega\times V. 
\end{equation}
	
{\bf Theorem 3.1.$\quad$}
\begin{theorem}  Let $\Omega\subset {\bf R^N}$ be a bounded domain. Let $u\in USC(\Omega\times V)$,  $w\in LSC(\Omega\times V)$ be respectively a sub and a super solution of (\ref{stationary}). Assume that (\ref{mass}), (\ref{b}), (\ref{g}) hold, and 
$$
	u\leq w \qquad \hbox{on} \quad \p(\Omega\times V). 
$$
Then, $u\leq w$ in $\Omega\times V$. \\
\end{theorem}

{\bf Proof.} We use the argument by contradiction. Assume that $\max_{\Omega\times V}(u-w)(x,v)$$=(u-w)(x_0,v_0)>0$ for $(x_0,v_0)$$\in \Omega\times V$. 
For $\a>0$,  put
$$
	\Phi_{\a}(x,y,v,v')=u(x,v)-w(y,v')-\a|x-y|^2 - \a|v-v'|^2. 
$$
Let $(\hx_{\a},\hy_{\a},\hv_{\a}, \hv'_{\a})$ be the maximum point of $\Phi_{\a}$ in $(\Omega\times V)^2$.   It is known (see \cite{users}) that 
\begin{equation}\label{error}
	(\hx_{\a},\hv_{\a}),\quad (\hy_{\a},\hv'_{\a})\to (x_0,v_0); \quad \a |\hx_{\a}-\hy_{\a}|^2, \quad \a |\hv_{\a}-\hv'_{\a}|^2 \to 0 \quad \hbox{as}\quad \a \to \infty. 
\end{equation}
 In the following, we abbreviate the index $\a$ for simplicity.
From the definition of the viscosity solution, for $p=2\a(\hx-\hy)$, 
$$
	\l u (\hx,\hv)+ \sup_{\a \in \mathcal{A}} \{ -\left\langle b(\hx,\hv,\a), p \right\rangle\} - 
	\int_{V} [u(\hx,v') - u(\hx,\hv)] dq(v') - g(\hx,\hv)\leq 0,
$$
$$
	\l w (\hy,\hv')+ \sup_{\a \in \mathcal{A}} \{ -\left\langle b(\hy,\hv',\a), p \right\rangle\} - 
	\int_{V} [w(\hy,v') - w(\hy,\hv')] dq(v') - g(\hy,\hv')\geq 0. 
$$
For any $\e>0$, we can take a control $\overline{\a}$ such that 
$$
	\l u (\hx,\hv) -\left\langle b(\hx,\hv,\overline{\a}), p \right\rangle  - 
	\int_{V} [u(\hx,v') - u(\hx,\hv)] dq(v') - g(\hx,\hv)\leq 0,
$$
$$
	\l w (\hy,\hv') -\left\langle b(\hy,\hv',\overline{\a}), p \right\rangle - 
	\int_{V} [w(\hy,v') - w(\hy,\hv')] dq(v') - g(\hy,\hv')\geq -\e. 
$$
By taking the difference of the above two inequalities, we get 
$$
	\l (u (\hx,\hv)-w (\hy,\hv')) \leq \e+ g(\hx,\hv) - g(\hy,\hv')\qquad\qquad\qquad\qquad\qquad\qquad
$$
$$
 	+ \left\langle b(\hx,\hv,\overline{\a}) - b(\hy,\hv',\overline{\a}), p \right\rangle 
	+ \int_{V} [u(\hx,v') - u(\hx,\hv) - w(\hy,v') +w(\hy,\hv') ] dq(v'). 
$$
 Since 
$$
	u (\hx,\hv)-w (\hy,\hv')- \a|\hx-\hy|^2 - \a|\hv-\hv'|^2 \qquad\qquad\qquad\qquad\qquad\qquad
$$
$$
	\geq u (\hx,v')-w (\hy,v')- \a|\hx-\hy|^2 - \a|v'-v'|^2 \quad \forall v'\in V, 
$$
by introducing this to the above,  from (\ref{b}), (\ref{g}), 
$$
	\l (u (\hx,\hv)-w (\hy,\hv')) \leq C (|\hx-\hy|^{\theta}+ |\hv-\hv'|^{\theta})
 	+ \a|\hx-\hy|(|\hx-\hy|+ |\hv-\hv'|)
$$
$$
	\leq C (|\hx-\hy|^{\theta}+ |\hv-\hv'|^{\theta}) + 2\a|\hx-\hy|^2+ \a |\hv-\hv'|^2. 
$$
From (\ref{error}), the right hands side of the above tends to zero as $\a$ goes to $\infty$. 
This contradicts to $\max_{\Omega\times {\bf R^N}}(u-w)(x,v)$$>0$, and we proved the claim. \\

{\bf Theorem 3.2.$\quad$}
\begin{theorem}  Assume that (\ref{mass}), (\ref{b}), and (\ref{g}) hold.
Then, there exists a unique solution $u(x,v)$ of (\ref{stationary})-(\ref{bc}). 
\end{theorem}

{\bf Proof.}  From the comparison result in Theorem 3.1, the existence of the solution is derived by the Perron's method (see \cite{users} for more details). 
Let $M=\max\{ |g|_{L^{\infty}(\Omega\times {\bf R^N})},   |\psi|_{L^{\infty}(\Omega\times {\bf R^N})}\}. $ Put $\uu(x,v)=-M$,  $\ou(x,v)=M$ for any $(x,v)$$\in \Omega\times {\bf R^N}$. It is clear that $\uu$, $\ou$ are respectively a sub and a super solution of (\ref{stationary})-(\ref{bc}). Let
$$
	\overline{U}(x,v)=\sup \{U(x,v)|\quad U\quad \hbox{is a subsolution of  (\ref{stationary})-(\ref{bc})}, \quad \uu\leq U\leq \ou\}. 
$$
Then, by the standard argument (\cite{users}), $\overline{U}(x,v)$ is a viscosity solution of  (\ref{stationary})-(\ref{bc}). The uniquness follows from  Theorem 3.1.\\

\end{document}